# An Alternative Distributed Control Architecture for Improvement in the Transient Response of DC Microgrids

Tuyen V. Vu, *Student Member, IEEE*, Sanaz Paran, *Student Member, IEEE*, Fernand Diaz-Franco, *Student Member, IEEE*, Touria El-Mezyani, *Member, IEEE* and Chris S. Edrington, *Senior Member, IEEE*

*Abstract*—Distributed secondary control plays an important role in DC microgrids, since it ensures system control objectives, which are power sharing and DC bus voltage stability. Previous studies have suggested using a control architecture that utilizes a parallel secondary bus voltage and current sharing compensation. However, the parallel controllers have a mutual impact on each other, which degrades the transient performance of the system. This paper reports on an alternative distributed secondary control architecture and controller design process, based on small signal analysis to alleviate the mutual effect of the current sharing and bus voltage compensation, and to improve the transient response of the system. Experimental results confirm the improved transient performance in the current sharing control and DC bus voltage stability utilizing the proposed control architecture.

## I. INTRODUCTION

IN recent years, the development of power electronics technology has made DC microgrid architectures and control a promising area for researchers. DC microgrids possess competitive control advantages over AC microgrids, including low transmission loss and simple control algorithms because there is no reactive power flow, frequency regulation, and synchronization [1]-[3].

In DC microgrids, the traditional control methodology is voltage droop control. Different types of droop control are found in the literature [4]-[6]. Droop control is widely used to maintain the proportional current sharing between distributed resources by reducing their output voltage, following a predefined droop characteristic. Thus, the enhancement of droop control in the current sharing attenuates the DC bus voltage stability [7]-[14]. An increment in droop parameters,

Manuscript received January 20, 2016; revised May 11, 2016 and June 27, 2016; accepted July 25, 2016. This work was supported by the Office of Naval Research via grant awards: N00014-10-1-09.

The authors are with the Florida State University Center for Advanced Power Systems, Tallahassee, FL 32310, USA (email: tvu@caps.fsu.edu; paran.sanaz@gmail.com; fediazfr@gmail.com; tmezyani@gmail.com; edrington@caps.fsu.edu).

for instance, results in increased accuracy of the current sharing, but an increased bus voltage drop. Comparatively, a reduction in droop parameters results in a decreased bus voltage drop, but inaccuracy in current sharing. Consequently, advanced methods based on the secondary control architecture [7], [13], [15]-[18], have been introduced for improving the system's performance. Ordinarily, the control structure is the combination of non-proportional current sharing compensation and deviated bus voltage compensation, which both generate the change in one input voltage reference of each power electronic converter. This scheme suggests accurate current sharing and bus voltage stability. However, it has a drawback of having a conflicting interest between the two control inputs for voltage and current compensation. The two control inputs have a mutual impact on each other in the system in order to achieve their own objective. Specifically, the adjusted voltage causes an unexpected transient in the current sharing and similarly in the bus voltage, which degrades the performance of the transient response of the current sharing and the DC bus voltage restoration. Therefore, a cascade control structure, which employs power sharing control as an inner loop and employs bus voltage control as an outer loop, is proposed to alleviate the aforementioned mutual control effect.

In addition to proposing the cascade control structure, system modeling to support control design is implemented in this paper. Literature review indicates that the process from system modeling to controller design, based on the stability criteria, have yet to be fully developed. Models for DC microgrids, for example, have been developed for stability analysis utilizing root locus in [19]-[21]; nevertheless, they mainly focus on stability analysis for the existing droop control and controller parameters instead of defining the secondary controller parameters based on the stability criteria. Hence, a model, which reflects the relationship between the input and output of the system to enhance the control system design and instantiation, is of fundamental importance. This paper addresses the control design subsequent to the modeling requirement.

Along with the improvement of the droop control, control architectures for microgrids are gaining a lot of attention. Centralized control architectures have been proposed in [22]-[26]. This type of architecture proposes synchronous information and central control. Consequently, due to the large and complex nature of microgrids, including: DC distribution, distributed variable load structure, and power components, the centralized technique becomes impractical because of: the possible need for real-time optimizations, dynamical changes of the system under plug and play operation, and reliable information exchange [27]. Therefore, a distributed power management system is required to ensure plug and play operation and to achieve common objectives of the systems, such as current sharing and bus voltage stability [28]-[32]. In this paper, a distributed power management control architecture is introduced to employ the proposed cascade control scheme, which fulfills the system's objectives.

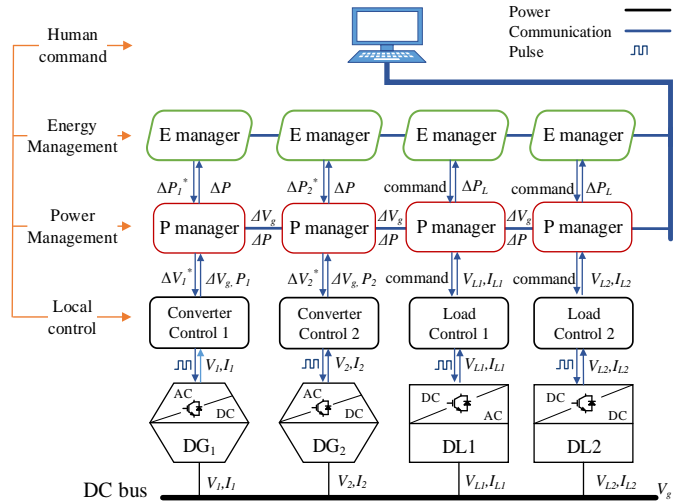

Fig. 1. Distributed power and energy management architecture.

This paper is organized as follows. Section II introduces a generalized distributed power and energy management architecture for a DC microgrid system, where the energy and power control tasks are classified hierarchically. The DC microgrid small signal model for two distributed power generations (DG) with power electronic converters is analyzed and formulated in Section III. The derived model leads to an alternative cascade distributed secondary control architecture, in which the power control is the inner loop, and the bus voltage control is the outer loop. Section IV presents a formal guideline for power and voltage control design in frequency domain using Bode plots. Section V details an experimental setup, in which the proposed control method is implemented and quantitatively compared against previous methods. Section VI summarizes the contributions of the paper.

## II. ARCHITECTURE OF DISTRIBUTED POWER MANAGEMENT

As mentioned, a distributed architecture, which employs the distributed control scheme is essential in the coordination of DG to maintain the operation of DC microgrids under various scenarios. Located in the distributed architecture are the distributed controllers, which implement a control scheme to ensure that the system is properly regulated. The distributed architecture and conventional management scheme will be discussed in the following sub-sections:

### A. Distributed Architecture

A control architecture that ensures the coordination of DG in microgrids can be found in a centralized manner [18]. Similar to the approach, but possessing a control and management flexibility, this paper proposes a distributed control architecture for a notional DC microgrid. The DC microgrid candidate presented in this paper includes: (1) two AC/DC rectifiers acting as DC power sources ($DG_1$ and $DG_2$); (2) a distributed AC load ($DL_1$) with an internal DC/AC inverter; and, finally, (3) a distributed DC load ($DL_2$) with an internal DC/DC converter. The DC microgrid with distributed control and management architecture is shown in Fig. 1.

In the distributed hierarchical architecture, a hierarchical controller is defined by three layers. Thus, each DG or load has three hierarchical controllers: a local controller, a power management controller, and an energy management controller. The distributed controllers in this architecture, called P-managers, are developed based on multi-agent technology [33]-[36]. The distributed power management control consequently is the distributed secondary control. The input to the P-managers is received from the energy management systems, appropriately called as E-manager. The E-managers determine the amount of energy supplied by each DG by means of generating the power command $\Delta P_i^*$ ($i = 1, 2$) to the P-manager (i.e. power agent) of the DG. The power command $\Delta P_i^*$ is the result of the energy management scheme applied in E-managers to achieve the system objective. Energy management schemes can be the optimization of the operational cost for renewable energy systems involving energy storage devices [37],[38], or fuel consumption minimization for operation of generators in ship power systems [39]. Since this paper focuses on the power management, the distributed energy management scheme is not further discussed. Although energy management scheme is not further discussed, to generate such a power command to the power manager $i$, a simple power reference calculation is selected based on the nominal power of the converters as

$$\Delta P_i^* = w_i \sum_{k=1}^{n=2} \Delta P_k \qquad (1)$$

where, the weight parameter $w_i$ is selected as

$$w_i = \frac{P_{ir}}{\sum_{i=1}^{n=2} P_{ir}} \qquad (2)$$

where, $P_{ir}$ is the rated power of $DG_i$.

In the distributed power management level, the P-managers, which are the focus of this paper, need to regulate the power and DC bus voltage simultaneously. The input of one P-manager $i$ in Fig. 1 is the power reference $\Delta P_i^*$ from the E-manager $i$; bus voltage deviation $\Delta V_{gi}$ and supplying power $\Delta P_i$ received from its local controller; and bus voltage deviation $\Delta V_g$ and power information $\Delta P_j$, ($j \neq i$) received from neighbors through a communication channel. The tasks of the P-manager are to follow the power command $\Delta P_i^*$ from the E-manager, and to minimize the DC bus voltage deviation $\Delta V_g$. The outputs of the P-manager are the voltage command $\Delta V_i^*$ sent to the local controller of the power electronic converter to perform the voltage and current regulation for the desired terminal voltage of the converter, which adjusts the power sharing and achieves the desired bus voltage stability.

### B. Conventional Power Management Scheme

To proper regulate the power sharing among DG and stabilize the DC bus voltage, a power management control scheme is required. The conventional power management control (secondary control) for a converter $i$ ($CONV_i$) shown in Fig. 2 varies the converter's voltage reference $V_i^{ref}$ as follows [7]:

$$V_i^{ref} = V_i^* - R_{di}I_i + \delta_{vi} + \delta_{Ii} \quad (3)$$

where $V_i^*$ is a fixed voltage reference of the converter, $I_i$ is the current sharing, $R_{di}$ is the virtual impedance, and $\delta_{vi}$ and $\delta_{Ii}$ are the outputs of the secondary voltage and current controllers. The virtual impedance $R_{di}$ is applied as the droop in the primary control. The secondary voltage and current controllers generate the voltage reference changes $\delta_v$, and $\delta_I$ in order to have the desired terminal voltage of the converter connecting to the DC bus.

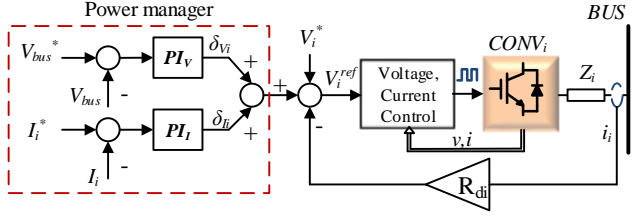

Fig. 2. Conventional power management scheme.

This control methodology in utilizing distributed P-manager utilizes additional outer bus voltage and current controllers $PI_V$ and $PI_I$. Indeed, the presence of two control inputs $\delta_v$ and $\delta_I$, which are the outputs of the voltage controller and the current controllers, have coupling effect on each other, and thus restrict the transient performance of the system. As a result, during the transient, the enhanced performance in the bus voltage restoration degrades the performance in the current sharing and vice versa. Hence, a necessary alternative distributed power management structure, based on small signal analysis that alleviates the conflict between the voltage and current compensation, is analyzed and proposed in the next Section.

## III. DC MICROGRID MODELING

To have a proper control design procedure to DC microgrids, modeling of these systems is required. DC microgrids modeling involves understanding the behavior of the system under critical disturbances, including input voltage variation and/or load change. As a result, a small signal model for a DC microgrid candidate is introduced for control design in this Section.

### A. Small Signal Model

Consider a voltage source $V$ connecting to a DC bus, the relationship of the converter in a microgrid can be seen as a voltage source interfacing with a constant voltage load $V_g$, see Fig. 3. Applying the Laplace transform to the Kirchhoff's circuit laws for the circuit in Fig. 3, the converter output voltage $V$ and power generation $P$ relationship is

$$P(s) = V_g \frac{V(s) - V_g}{R + Ls} \quad (4)$$

where $R$ and $L$ are the cable resistance and inductance, respectively. Consider a variation in converter output voltage $\Delta V$. This results in a change in the power $\Delta P$ generated to the DC microgrid as shown in (5). Thus, the transfer function (6) derived from (5) represents the relationship of the change of the output voltage to the change of the output power.

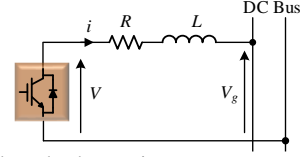

Fig. 3. Constant voltage load operation.

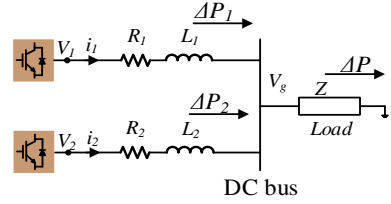

Fig. 4. Two DC sources microgrid.

$$P(s) + \Delta P(s) = V_g \frac{V(s) + \Delta V(s) - V_g}{R + Ls} \quad (5)$$

$$\frac{\Delta P(s)}{\Delta V(s)} = V_g \frac{1}{R + Ls} \quad (6)$$

The previous assumption of one converter connected to a constant voltage load is not adequate since microgrids are multi-terminal connected systems, coupling voltage and current. Thus, to formulate the power coupling consider a system consisting of two converters $V_1$ and $V_2$ supplying power to a generic load with impedance $Z$ as seen in Fig. 4. $R_1$, $R_2$ and $L_1$, $L_2$ are the cable resistances and inductances, respectively. Applying the Laplace transform to Kirchhoff's law for the above circuit, one obtains:

$$\frac{V_1(s) - V_g(s)}{R_1 + L_1 s} + \frac{V_2(s) - V_g(s)}{R_2 + L_2 s} = \frac{V_g(s)}{Z(s)} \quad (7)$$

where $Z(s)$ is the load impedance. Development of a small signal model, which supports the control system design, requires, the utilization of differentials for voltage $\Delta V_1$, $\Delta V_2$ and load disturbance $\Delta Z$; and their influence in the system including the bus voltage change $\Delta V_g$, and power sharing $\Delta P_1$, $\Delta P_2$. Thus, two assumptions are made for these parameter variations in the system (7).

Assumption 1: Small variations of the converter's output voltages $\Delta V_1$, $\Delta V_2$ occur while maintaining the load variation ($\Delta Z = 0$).

Assumption 2: Small variations in the load, $\Delta Z$ occur while maintaining the converter's output voltages ($\Delta V_1 = \Delta V_2 = 0$).

Assumption 1 results in the change of the bus voltage $\Delta V_g$. Substituting the changes into (7) yields (8).

$$\frac{V_1(s) + \Delta V_1(s) - V_g(s) - \Delta V_g(s)}{R_1 + L_1 s} + \frac{V_2(s) + \Delta V_2(s) - V_g(s) - \Delta V_g(s)}{R_2 + L_2 s} = \frac{V_g(s) + \Delta V_g(s)}{Z(s)} \quad (8)$$

Based on (7) and (8), the voltage change $\Delta V_g$ in the DC bus depends on $\Delta V_1$, $\Delta V_2$ via the relationship in (9).

$$\Delta V_g(s) = \frac{(R_2 + L_2 s)\Delta V_1(s) + (R_1 + L_1 s)\Delta V_2(s)}{R_1 + R_2 + (L_1 + L_2)s} \quad (9)$$

Assumption 2 results in the change of bus voltage $\Delta V_g$. Substituting the changes into (7) yields (10).

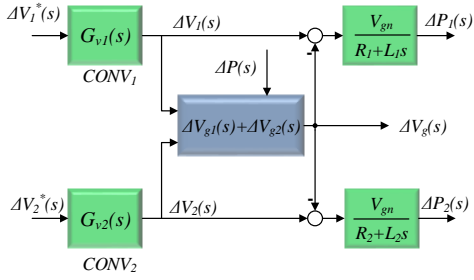

Fig. 5. DC microgrid coupling model.

$$\frac{V_1(s) - V_g(s) - \Delta V_g(s)}{R_1 + L_1 s} +$$

$$\frac{V_2(s) - V_g(s) - \Delta V_g(s)}{R_2 + L_2 s} = \frac{V_g(s) + \Delta V_g(s)}{Z(s) + \Delta Z(s)} \quad (10)$$

Based on (7) and (10) and multiplying by nominal grid voltage $V_{gn}$, (11) is derived as

$$\left(\frac{-\Delta V_g(s)}{R_1 + L_1 s} + \frac{-\Delta V_g(s)}{R_2 + L_2 s}\right) V_{gn}$$

$$= \left(\frac{V_g(s) + \Delta V_g(s)}{Z(s) + \Delta Z(s)} - \frac{V_g(s)}{Z(s) + \Delta Z(s)}\right) V_{gn} \quad (11)$$

It can be seen that the right side of (11) approximately represents the power change in the load $\Delta P$.

$$\left(\frac{-\Delta V_g(s)}{R_1 + L_1 s} + \frac{-\Delta V_g(s)}{R_2 + L_2 s}\right) V_{gn} = \Delta P \quad (12)$$

Therefore, the bus voltage drop is given by (13):

$$\Delta V_g(s) = -\frac{1}{V_{gn}} \frac{(R_1 + L_1 s)(R_2 + L_2 s)}{R_1 + R_2 + (L_1 + L_2)s} \Delta P(s) \quad (13)$$

As seen, there are two factors contributing to the bus voltage variation, which are the change in converter's output voltage and the change in the load. Therefore, in this linear analysis the superposition principle is utilized to determine the total bus voltage variation $\Delta V_g$ by combining voltage change caused by the change in converter's output voltage in (9) $\Delta V_g = \Delta V_{g1}$ and voltage change caused by load change in (13) $\Delta V_g = \Delta V_{g2}$. Thus, the bus voltage drop $\Delta V_g$ is expressed as

$$\Delta V_g(s) = \frac{(R_2 + L_2 s)\Delta V_1(s) + (R_1 + L_1 s)\Delta V_2(s)}{R_1 + R_2 + (L_1 + L_2)s} -$$

$$\frac{1}{V_{gn}} \frac{(R_1 + L_1 s)(R_2 + L_2 s)}{R_1 + R_2 + (L_1 + L_2)s} \Delta P(s) \quad (14)$$

The variations in power exchange between converters and loads are shown in (15).

$$\Delta P_1(s) = \frac{\Delta V_1(s) - \Delta V_g(s)}{R_1 + L_1 s} V_{gn}$$

$$\Delta P_2(s) = \frac{\Delta V_2(s) - \Delta V_g(s)}{R_2 + L_2 s} V_{gn} \quad (15)$$

It is noted that the derived small signal models, which are the relationship between the bus voltage drop $\Delta V_g$ and the power output of the converters $\Delta P$, and the relationship between the power sharing $\Delta P_1$, $\Delta P_2$ and output voltages of the converters $\Delta V_1$, $\Delta V_2$, will subsequently be utilized for the proposed control system design. The derived small signal model for the DC microgrid based on the bus voltage drop $\Delta V_g$ and power exchange $\Delta P_1$ and $\Delta P_2$ is graphically illustrated via the block diagram with the transfer functions as shown in

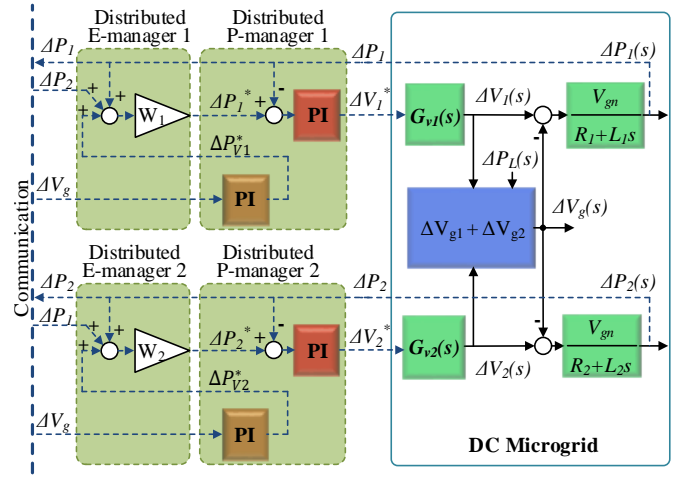

Fig. 6. Distributed power control diagram.

Fig. 5, where $\Delta V_1^*$ and $\Delta V_2^*$ are changes in voltage references and $G_{V1}(s)$ and $G_{V2}(s)$ are the voltage transfer functions for converter 1 ($CONV_1$) and converter 2 ($CONV_2$), respectively.

### B. Proposed Distributed Power Management

The relationship in (14) between the bus voltage drop $\Delta V_g$ and load change $\Delta P$ indicates that the generated power by the converters is taken into account for the DC bus voltage restoration. In addition, the relationship in (15) shows that the converter output voltage is the input for the power sharing regulation. As such, the power management diagram in the hierarchical distributed scheme with communication for the microgrid is proposed as shown in Fig. 6. In this scheme, the local controllers regulate the converter current and voltage control. There is no droop control in the primary control level.

The distributed P-manager, after receiving the power command from E-manager, regulates power sharing and bus voltage in the microgrid. The P-manager is the cascade control of the inner power control loop and the outer DC bus voltage loop. The output of the DC bus voltage PI controller is the power required $\Delta P_{Vi}^*$ to modify the voltage of the converters, which is summed up with the power demand in order to modify the power reference equation for each converter in (1) in the microgrid. The new power reference calculation $\Delta P_i^*$ for the converter $i$ ($i = 1, 2$) is defined as follows:

$$\Delta P_i^* = w_i \left(\sum_{i=1}^{n=2} \Delta P_i + \Delta P_{Vi}^*\right) \quad (16)$$

where $\Delta P_i$ is the power feedback from the sending end of converter $i$, $\Delta P_{Vi}^*$ is the power reference generated by the bus voltage controller, and $w_i$ is the weight parameter selected by (2) and it is implemented in the distributed E-manager $i$. The output of each of the power PI controllers are the reference voltages $\Delta V_i^*$ ($i = 1, 2$) for local control of each converter.

## IV. MANAGEMENT CONTROLLER DESIGN GUIDELINES

### A. Power Controller

The open loop small signal relationships between the input voltage and output power of each converter derived from block diagram in Fig. 6 are defined as

$$G_{OP1}(s) = \frac{\Delta P_1(s)}{\Delta V_1^*(s)} = \frac{\Delta P_1(s)}{\Delta V_1(s)} \frac{\Delta V_1(s)}{\Delta V_1^*(s)} = G_{v1}(s) \frac{V_{gn}}{R_1 + L_1 s}$$

$$G_{OP2}(s) = \frac{\Delta P_2(s)}{\Delta V_2^*(s)} = \frac{\Delta P_2(s)}{\Delta V_2(s)} \frac{\Delta V_2(s)}{\Delta V_2^*(s)} = G_{v2}(s) \frac{V_{gn}}{R_2 + L_2 s} \quad (17)$$

where $G_{vi}$ ($i = 1,2$) is the voltage-loop transfer function of the converter $i$. Deriving the system-level models for microgrids, fast dynamics of power converters with multiple-order in the transfer function $G_{vi}$ are neglected. Consequently, $G_{vi}$ can be considered as a delay, which is equivalent to the following reduced first-order model:

$$G_{vi} = \frac{1}{1 + \tau_i s} \quad (18)$$

where $\tau_i$ represents the time delay of the voltage control loop. To regulate the power flowing from converter $i$ to the DC bus, a PI controller is utilized, see Fig. 7. Based on the open loop transfer function between reference voltage of converter $\Delta V_i^*$ with power output $\Delta P_i^*$, ($i = 1,2$) (17), the analysis in the frequency domain is utilized for the PI controller's design based on the phase margin and crossover frequency requirements of the control system. The detailed procedure is conducted and validated in the Section V.

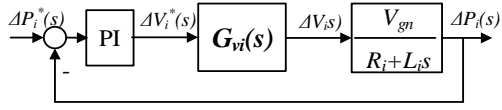

Fig. 7. Power control loop for one energy resource.

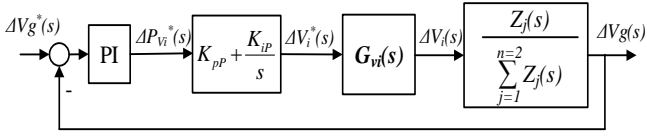

Fig. 8. Bus voltage control loop for one energy resource.

### B. Bus Voltage Controller

Based on the input and output relationships between signals in Fig. 6, small signal models derived from the block diagram for the bus voltage control between the input power reference $\Delta P_{Vi}^*$ and output bus voltage $\Delta V_g$ of each converter are

$$G_{OVg1}(s) = \frac{\Delta V_g(s)}{\Delta P_{V1}^*(s)} = (K_{pP} + \frac{K_{iP}}{s})G_{v1}(s) \frac{Z_2(s)}{Z_1(s) + Z_2(s)}$$

$$G_{OVg2}(s) = \frac{\Delta V_g(s)}{\Delta P_{V2}^*(s)} = (K_{pP} + \frac{K_{iP}}{s})G_{v2}(s) \frac{Z_1(s)}{Z_1(s) + Z_2(s)} \quad (19)$$

where $\Delta P_{V1}^*$ and $\Delta P_{V2}^*$ are the power changes, resulting in the bus voltage deviation $\Delta V_g$; $G_{vi}(s)$ is the closed-loop transfer function of the local voltage control loop of $CONV_i$ (18); $K_{pP}$ and $K_{iP}$ are the power PI controller parameters; and $Z_1$ and $Z_2$ are the cable impedances. The control diagram for converter $i$ of with PI controller is shown in Fig. 8. Similar to the power control design in the previous part, based on the small signal relationship between the power reference $\Delta P_{Vi}^*$ and the bus voltage deviation $\Delta V_g$, an analysis in the frequency domain utilizing Bode plot is implemented based on the phase margin and crossover frequency design criteria of the bus voltage restoration control system. Details about the design procedure for a study case are explained in Section V.

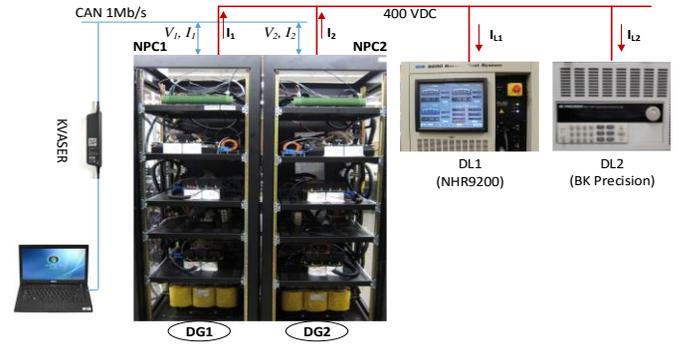

Fig. 9. Experimental setup.

## V. CASE STUDY

To compare and validate the effectiveness of the proposed control algorithm with control design procedure, an experimental system that represents a DC microgrid is setup. Then, the power management controllers are designed. Next, the control parameters are then tested to verify the transient performance of the proposed algorithm. Subsequently, to compare and contrast the proposed control against the conventional control, the conventional control is implemented using the same parameters and operating conditions; where the voltage controller's gains are varied within a range of applicable values. The purpose of the experimental study is to verify that there is a mutual effect between the current and voltage controllers in the conventional scheme, which is mitigated in the proposed controller.

### A. System Description

The system illustrated in Fig. 9 is a 400V DC microgrid, which includes two Neutral Point Clamped (NPC) converters that are controlled in order to share the power proportional to their rated power ($P_{NPC1}:P_{NPC2} = 2:1$). The two NPC are powered by two AC transformers, which connect to the same AC source in the laboratory. In details, 4kW is assumed rated power for $NPC_1$, and 2kW is assumed rated power for $NPC_2$. The distributed loads connected to the DC bus are 4kW (NHR9200) and 2kW (BK Precision DC load), respectively. The communication between two DSPs, TMS28335, for each $NPC$ is achieved via CAN at a rate of 1Mb/s. The information exchanged between the DSPs are the terminal bus voltages $V_{NPC1}$ and $V_{NPC2}$, and the output currents $I_{NPC1}$ and $I_{NPC2}$. The data acquisition and control system activation commands via CAN using the computer is implemented through the Kvaser Leaf Light v2. The three-level active rectifier (NPC) with the control algorithm instantiated in the d-q frame [40] is shown in Fig. 10. The system parameters for the DC microgrid are listed in TABLE I.

### B. Distributed Controllers Design

#### 1) Power controller design

The open-loop relationship in (17) between the change in the output power $\Delta P_i$ and the change in the input voltage reference change $\Delta V_i^*$ ($i = 1,2$), has a crossover frequency of $\omega_{0c} = 116.6 \ rad/s$. The criterion is to have the system response as approximately fast as the open-loop system in the closed-loop design. Thus, the crossover frequency for the controller design is selected as $\omega_{Pc} = 100 \ rad/s$. The phase margin is chosen as $\varphi = 70^0$ to ensure that the closed-loop system is stable under disturbances or uncertainties. Based on the frequency response (Fig. 11a) with the specified criteria, the desired power PI controllers' parameters for the NPC are $K_{pP} = 0.001$ and $K_{iP} = 0.130$.

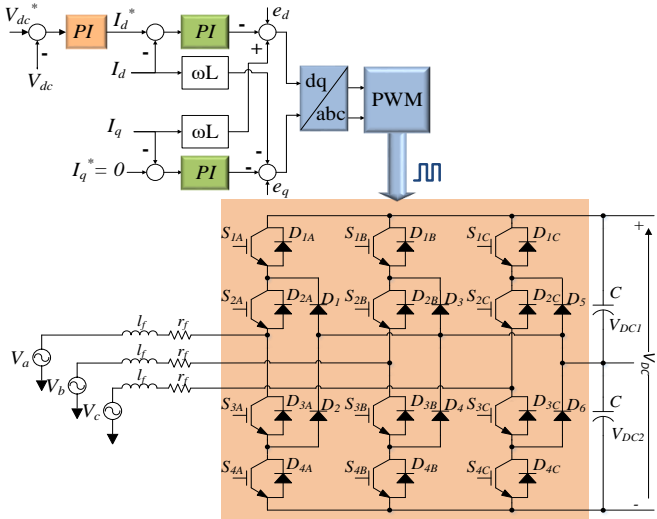

Fig. 10. NPC topology and control diagram.

TABLE I
SYSTEM PARAMETERS

| Symbol | Quantity | Values |
|--------|----------|--------|
| $V_{L-L}$ | NPC Input voltage | 208 V (60Hz) |
| $V_{DC}$ | NPC Output voltage | 400 V |
| $l_i$ | NPC input filter inductor | 2.07 mH |
| $f$ | Switching frequency | 20 kHz |
| $C$ | NPC output capacitor | 380 µF |
| $r_{fi}$ | NPC input resistor | 0.2 Ω |
| $R_i$ | Cable resistance | 0.5 Ω |
| $L_i$ | Cable inductance | 3 mH |
| $\tau_i$ | NPC voltage loop time constants | 0.005 s |

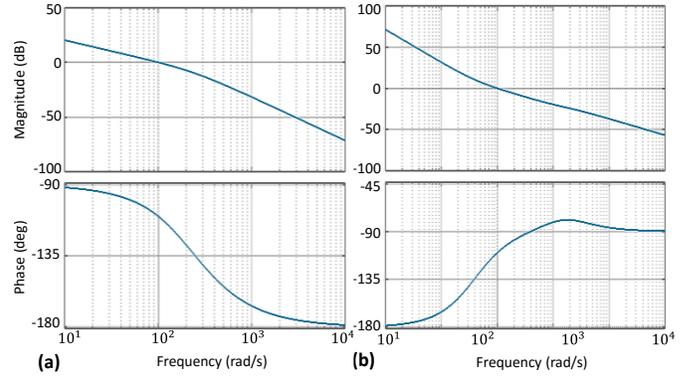

Fig. 11. Bode plots of one converter. (a) power control, (b) bus voltage control.

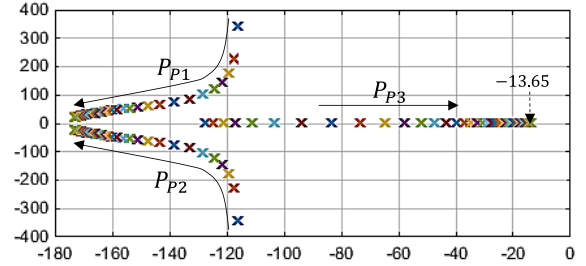

Fig. 12. Root locus of power control (Increasing in cable impedance).

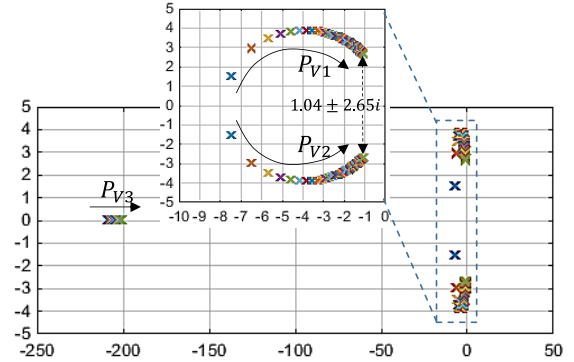

Fig. 13. Root locus of bus voltage control (Increasing in cable impedance).

### 2) Bus voltage controller design

The bus voltage control, which is the outer loop of the power control loop, has a slower transient in comparison to the power control. In cascade control, the bandwidth of the outer loop is selected as 10 times approximately smaller than the one of the inner loop. Thus, the bus voltage control crossover frequency is selected as $\omega_{Vc} = 0.1\omega_{pc} = 10 \ rad/s$. The phase margin is chosen as $\varphi = 70^0$ to ensure that the bus voltage PI controllers' for the NPC are derived as $K_{pV} = 142.9$ and $K_{iV} = 563.8$.

### C. Stability Analysis

There are possible changes in system structure, which result in system's parameters changes (cable impedance changes). These changes possibly destabilize the system with designed control parameters. Thus, it is necessary to verify the benefit of the proposed technique in cases, where there are variations in the cable impedances of the system. To analyze the stability effect, root-locus is utilized for developed small signal models in (17) and (19) with their designed control parameters. In this analysis, cable impedance variation in one of the two voltage sources (first voltage source) is taken into account.

The case assumes that there are changes in cable impedance, but the $R_1/L_1$ ratio is assumed to be fixed as $R_1/L_1 = 0.5/0.003$. Consider that the maximum voltage regulation ratio between the sending end at the source terminal and the receiving end at the load bus terminal is 5% at 10A rated current supplying from the source. Therefore, the maximum resistance $R_{1max}$ can be changed as $R_{1max} = 0.05 \times V_{gn}/10 = 2\Omega$. Suppose that different scenarios result in the cable resistance changes in the first converter as $R_1$ varies between $0.1\Omega$ and $2\Omega$. These changes result in the change of cable inductance $L_1$ as it varies from 0.6mH to 12mH because the $R_1/L_1$ ratio is assumed to be fixed. Applied these changes to plot the root locus of the models shown in (17) and (19).

As seen in Fig. 12, as the impedance increases, the pole $P_{P3}$ moves toward the imaginary axis, and it becomes more dominant than the poles $P_{P1}$ and $P_{P2}$. The movement of $P_{P3}$, which is terminated at the value of -13.65 illustrated that the designed control ensures the stability of the power control loop in the range of impedance changes. In bus voltage regulation, Fig. 13 shows that the two dominant poles $P_{V1}$ and $P_{V2}$ move toward the imaginary axis but they terminate at $P_{V1} = 1.04 + 2.65i$ and $P_{V2} = 1.04 - 2.65i$. The poles movement within the left half plane of the imaginary axis demonstrated that the bus voltage control is stable. Thus, the proposed method is stable and robust to the cable impedance changes or the changes in network structure of the system.

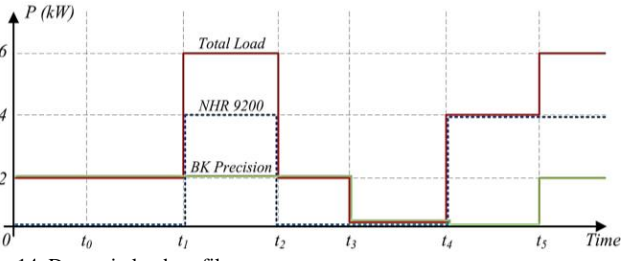

Fig. 14. Dynamic load profile.

### D. Experimental Results and Analysis

The expected results of the experiments are (1) 2:1 current sharing ratio between $NPC_1$ and $NPC_2$, and (2) 400 V bus voltage operation. The efficacy of the proposed control methodology is demonstrated via test cases, which all utilize the constant-power load profile, shown in Fig. 14. The solid red line represents the total load while the blue dash and solid yellow lines represent two distributed loads (NHR9200 and BK Precision).

The test cases conducted for comparisons between the conventional method and the proposed method are as follow: The first case utilizes the proposed control method, and the results for the terminal voltages $V_{NPC1}, V_{NPC2}$, and the current sharing between converters $I_{NPC1}, I_{NPC2}$ are illustrated in Fig. 15. To observe and analyze the mutual effect of the voltage and current controllers under the conventional control method, two other test cases utilizing conventional method are constructed. The first test case employed the conventional control architecture with the low gains in the bus voltage controller, which are selected as $K_{pV} = 0.2$, and $K_{iV} = 1$; to verify that the current sharing control is not affected by the bus voltage controllers with low gains. The results of this case are presented in Fig. 16. The second test case makes use of the conventional control scheme with high voltage controller gains, which are chosen as $K_{pV} = 1$, and $K_{iV} = 20$; and results are shown in Fig. 17. This test case illustrates the improvement in the voltage control while simultaneously demonstrating the degradation in the current sharing in the conventional control's performance. Note that in all the scope plots, which were recorded using a Yokagama DL850 that there are 100 V/major division and 5 A/major division.

The experimental data is analyzed for the comparison between the proposed method and the conventional method in order to investigate the system's behavior during the power management control activation, and load increment. Since there is high frequency noise in current measurement probe, the current data is passed through a low pass filter in Matlab to analyzing the transient improvement of the proposed method.

The results from proposed method and conventional method with low gains in the bus voltage controllers are compared through Fig. 18 and Fig. 19. Current sharing results are shown in Fig. 18a and Fig. 19a, and terminal bus voltage restoration is shown in Fig. 18b and Fig. 19b, respectively. The notations $I_{NPCi}new$, $V_{NPCi}new$, $I_{NPCi}old$, and $I_{NPCi}old$, $i = (1, 2)$ represent current supply and bus voltage of converters in the proposed method and conventional method, respectively. These figures indicate that the low gains in the bus voltage controller in the conventional method gives a comparable current response with the proposed method, because the low gains in the voltage control has only a small impact on the current control. However, the low gains of bus voltage controllers also result in a slower bus voltage transient response in comparison to the proposed method. Specifically, in the case of the power management control activation at $t = 5s$, the proposed method takes 0.25s to reach to nominal voltage value, while the conventional method takes 3.0s to reach to the nominal voltage value. In the event of a total load increment from 0 kW to 4kW at $t = 20s$, the proposed method takes 0.4s to restore the bus voltage, while the conventional method takes 3s to do the same task.

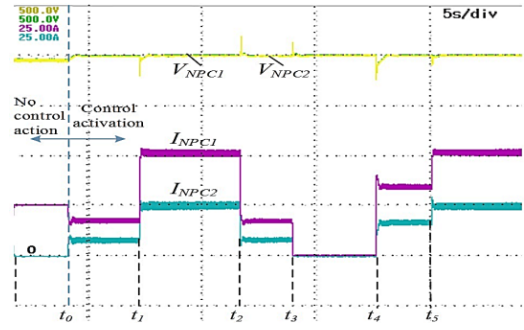

Fig. 15. Current and bus voltage profile of the proposed control algorithm.

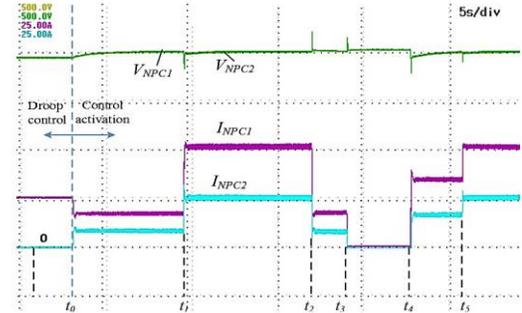

Fig. 16. Current and voltage profile with low gains in the voltage controller of the conventional algorithm.

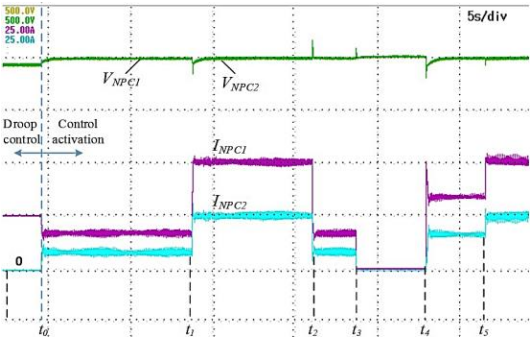

Fig. 17. Current and voltage profile with high gains in the voltage controller of the conventional algorithm.

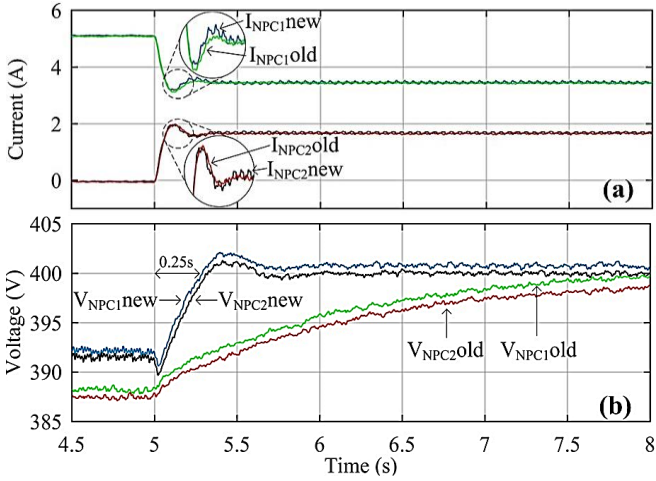

Fig. 18. Transients in the proposed method, and conventional method with low gains in the bus voltage controller during the control activation. (a) current sharing, (b) bus voltage.

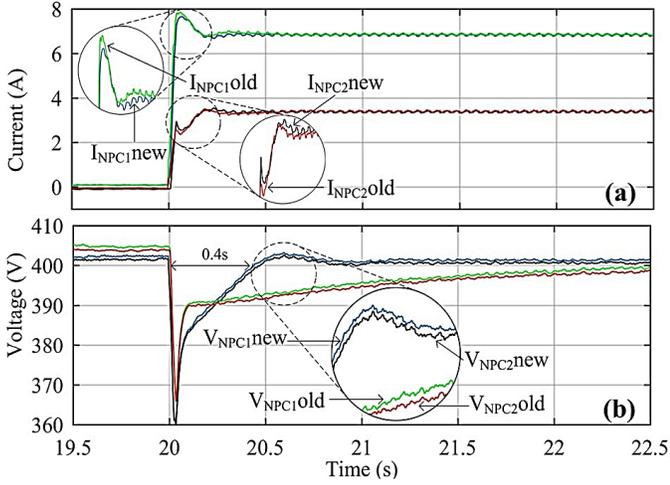

Fig. 19. Transients in the proposed method and conventional method with low gains in the bus voltage controller during a load increment. (a) current sharing, (b) bus voltage.

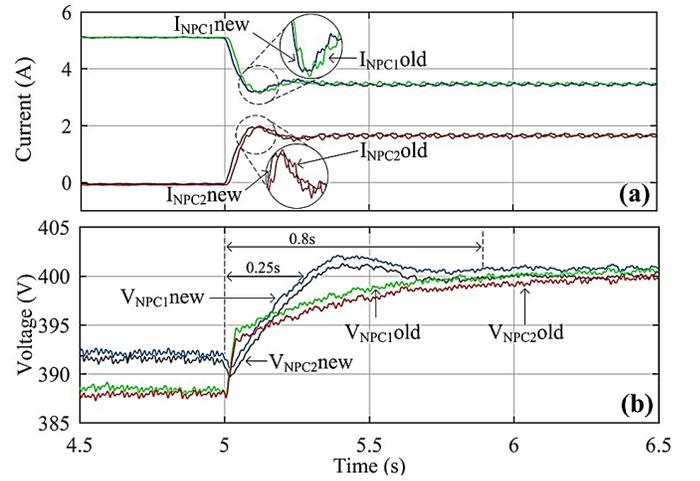

Fig. 20. Transients in the proposed method and conventional method with high gains in bus voltage controller during the control activation. (a) current sharing, (b) bus voltage.

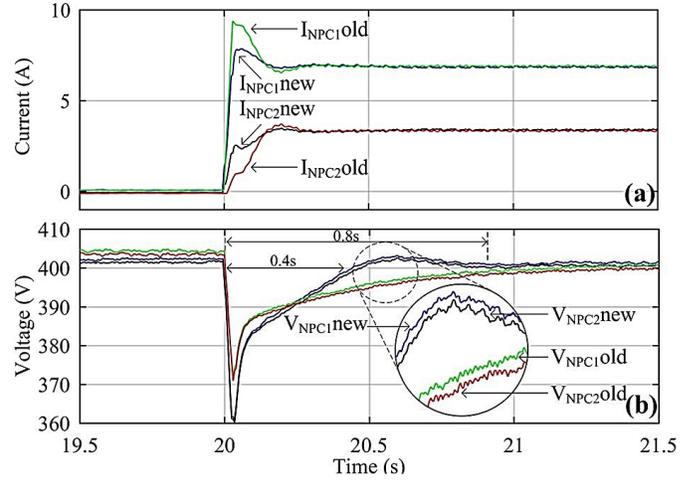

Fig. 21. Transient in the proposed method and conventional method with high gains in the bus voltage controller during a load increment. (a) current sharing, (b) bus voltage.

As seen, the low gains in the bus voltage control of the conventional method have less effect in the current sharing control but cause a slower bus voltage restoration. To improve the bus voltage control, high gains in bus voltage control are implemented. Fig. 20 and Fig. 21 present the comparisons between the proposed method and the conventional method with high gains in the bus voltage controller. The current sharing results shown in Fig. 20a, and Fig. 21a; and bus voltage control shown in Fig. 20b, and Fig. 21b, respectively, represent the comparison between the conventional control method with high gains in the DC bus voltage controllers and the proposed method. The same events are implemented as the conventional and the proposed control are activated at 5s; and the total load change from 0kW to 4kW is instantiated at 20s. In both events, the conventional control with high gains in the DC bus voltage controllers achieves a reduced 0.8s settling time compared to 3.0s settling time of bus voltage controllers with the previous low gains, improving the bus voltage transient response. However, the 0.8s settling time achieved by increasing the gains in the bus voltage controller is larger than the 0.4s settling time obtained by the proposed method.

Additionally, a higher overshoot for $NPC_1$ is observed in the current transient during the load increment event (Fig. 21a). Therefore, the improvement in the bus voltage controller of the conventional method by increasing the gains also simultaneously degrades the performance of the current sharing control. Comparatively, the proposed method enhances the current sharing control.

To quantify the effectiveness of the proposed method, the integral of time and absolute error (ITAE) criterion is utilized for bus voltage and current sharing control assessment. The bus voltage control and current sharing control are assessed by $ITAE_V$ (20), and $ITAE_I$ (21) respectively.

$$ITAE_V = \int_0^T t\left(V_{ref} - \frac{V_{NPC1} + V_{NPC2}}{2}\right) dt \quad (20)$$

$$ITAE_I = \int_0^T t\left[(I_{ref1} - I_{NPC1}) + (I_{ref2} - I_{NPC2})\right] dt \quad (21)$$

where T is the evaluation time selected as $T = 2s$; $V_{ref} = 400V$ is the bus voltage reference; $I_{ref1}$ and $I_{ref2}$ are the current sharing references; and $t$ is time. The ITAE comparison in TABLE II illustrates that there is a trade-off between the current and voltage control in the conventional method, in which the enhanced voltage response decreases the performance of the current sharing response and in contrast, that an improved current response degrades the bus voltage response. The improvement of the proposed methodology is illustrated via the reduction in the ITAE.

TABLE II
INTEGRAL OF TIME AND ABSOLUTE ERROR COMPARISON

| Case | | Conventional method (low-gain controllers) | Conventional method (high-gain controllers) | Proposed method |
|------|------|------|------|------|
| Control action (at t = 5s) | $ITAE_V$ | 1651 | 248.192 | 193.572 |
| | $ITAE_I$ | 232.969 | 304.371 | 225.574 |
| Load increment (at t = 20s) | $ITAE_V$ | 1594 | 697.358 | 563.225 |
| | $ITAE_I$ | 475.668 | 534.113 | 420.472 |

Consequently, the qualitative assessments through Fig. 18-Fig. 21 and the quantitative analysis using ITAE shown in TABLE II demonstrate that the proposed control method improves both current sharing and bus voltage stability of the DC microgrid simultaneously.

## VI. CONCLUSION

This paper addresses the transient response of the conventional power management control from the viewpoint of the tradeoff between secondary voltage and current controller performance. The response is enhanced by employing an alternative power management methodology derived from the small signal model of the DC microgrid. The controllers' design procedure, based on the analysis in the frequency domain, is conducted for the specified phase margins and crossover frequencies of the power and bus voltage control systems. An extended stability study for variation in distribution cable impedances is conducted to verify the robustness of the proposed control algorithm. The experimental data and comparative study confirm the improvement in transient performance in the current sharing and bus voltage stability simultaneously when utilizing the proposed power management control.

## REFERENCES

[1] J. J. Justo, F. Mwasilu, J. Lee, and J. W. Jung, "AC-microgrids versus DC-microgrids with distributed energy resources: A review," *Renewable and Sustainable Energy Reviews*, vol. 24, pp. 387–405, 2013.

[2] R. S. Balog, W. W. Weaver, and P. T. Krein, "The load as an energy asset in a distributed DC smartgrid architecture," *IEEE Trans. Smart Grid*, vol. 3, no. 1, pp. 253–260, Feb. 2012.

[3] A. Kwasinski, "Quantitative evaluation of DC Microgrids Availability: Effect of System Architecture and Converter Topology Design Choices," *IEEE Trans. Power Electron.*, vol.26, no.3, pp. 835 – 851, March 2011.

[4] S. Luo, Z. Ye, R. L. Lin, and F. C. Lee. "A classification and evaluation of paralleling methods for power supply modules," In *Proceedings of Power Electronics Specialists Conference*, vol. 2, pp. 901–908, 1999.

[5] L. Batarseh, K. Sin, and H. Lee, "Investigation of the Output Droop Characteristics of Parallel-Connected DC-DC Converters," *Power Electronics Specialists Conference*, vol. 2, pp. 1342-1351, Jun. 1994.

[6] J. Perkinson, "Current Sharing of Redundant DC-DC Converters in High Availability Systems-A Simple Approach," *Applied Power Electronics Conference and Exposition*, APEC '95, pp. 952-956, 1995.

[7] X. Lu, J. M. Guerrero, K. Sun, and J. C. Vasquez, "An Improved Droop Control Method for DC Microgrids Based on Low Bandwidth Communication With DC Bus Voltage Restoration and Enhanced Current Sharing Accuracy," *IEEE Trans. Power Electron.*, vol. 29, no. 4, pp. 1800–1812, April 2014.

[8] T. V. Vu, "Distributed Robust Adaptive Droop Control for DC Microgrids," Ph.D. dissertation, Dept. Elect. Eng., Florida State Univ., Tallahassee, FL, USA, 2016.

[9] T. Dragicevic, J. M. Guerrero, and J. C. Vasquez, "A Distributed Control Strategy for Coordination of an Autonomous LVDC Microgrid Based on Power-Line Signaling," *IEEE Trans. Ind. Electron.*, vol. 61, no. 7, pp. 3313–3326, July 2014.

[10] X. Lu, J. M. Guerrero, K. Sun, J. C. Vasquez, R. Teodorescu, and L. Huang, "Hierarchical Control of Parallel AC-DC Converter Interfaces for Hybrid Microgrids," *IEEE Trans. Smart Grid*, vol. 5, no. 2, pp. 683–692, March 2014.

[11] M. Hamzeh, A. Ghazanfari, Y.A.-R.I. Mohamed, "Modeling and Design of an Oscillatory Current-Sharing Control Strategy in DC Microgrids," *IEEE Trans. Ind. Electron.*, vol. 62, no. 11, pp. 6647 - 6657, Nov. 2015.

[12] C. Jin, P. Wang, J. Xiao, Y. Tang, and F. H. Choo, "Implementation of Hierarchical Control in DC Microgrids," *IEEE Trans. Ind. Electron.*, vol. 61, no. 8, pp. 4032–4042, Aug. 2014.

[13] S. Anand, B. G. Fernandes and J. M. Guerrero, "Distributed Control to Ensure Proportional Load Sharing and Improve Voltage Regulation in Low Voltage DC Microgrids," *IEEE Trans. Power Electron.*, vol. 28, no. 4, pp. 1900–1913, April 2013.

[14] K. Rouzbehi, A. Miranian, J. I Candela, A. Luna, and P. Rodriguez, "A Generalized Voltage Droop Strategy for Control of Multiterminal DC Grids," *IEEE Trans. Industry Applications*, vol. 51, no. 1, pp. 607–618, Feb. 2015.

[15] V. Nasirian, S. Moayedi, A. Davoudi, and F. Lewis, "Distributed cooperative control of DC microgrids," *IEEE Trans. Power Electron.*, vol. 30, pp. 2288–2303, April 2015.

[16] Q. Shafiee, J. M. Guerrero, and J. Vasquez, "Distributed secondary control for islanded microgrids—A novel approach," *IEEE Trans. Power Electron.*, vol. 29, no. 2, pp. 1018–1031, Feb. 2014.

[17] V. Nasirian, A. Davoudi, and F. L. Lewis, "Distributed adaptive droop control for DC microgrids," in *Proc. IEEE Appl. Power Electron. Conf. Expo.*, pp. 1147–1152, Mar. 2014.

[18] J. M. Guerrero, J. C. Vasquez, J. Matas, L. G. de Vicuna, and M. Castilla, "Hierarchical control of droop-controlled AC and DC microgrids–a general approach toward standardization," *IEEE Trans. Ind. Electron.*, vol. 58, no. 1, pp. 158–172, Jan. 2011.

[19] S. Anand and B. G. Fernandes, "Steady state performance analysis for load sharing in DC distributed generation system," in *Proc. 10th Int. Conf. Environ. Electr. Eng.*, pp. 1–4, May 2011.

[20] S. Anand, and B. Fernandes, "Reduced order model and stability analysis of low voltage DC microgrid," *IEEE Trans. Ind. Electron.*, vol. 60, no. 11, pp. 5040–5049, Nov. 2013.

[21] Q. Shafiee, T. Dragicevic, J. C. Vasquez, and J. M. Guerrero, "Modeling, stability analysis and active stabilization of multiple DC-microgrid clusters," *IEEE International Energy Conference*, May 2014.

[22] J. A. Peças Lopes, C. L. Moreira, and A. G. Madureira, "Defining control strategies for analysing microgrids islanded operation," in *Proc. IEEE Power Technology*, Russia, pp. 1–7, June 2005.

[23] A. Tsikalakis and N. Hatzargyriou "Centralized control for optimizing microgrids operation," *IEEE Trans. Energy Convers.*, vol. 23, no. 1, pp. 241–248, Mar. 2008.

[24] A. G. Tsikalakis and N. D. Hatziargyriou, "Centralized control for optimizing microgrids operation," *IEEE Trans. Energy Convers.*, vol. 23, no. 1, pp. 241–248, Mar. 2008.

[25] T. V. Vu, S. Paran, F. Diaz, T. E. Meyzani and C. S. Edrington, "Model predictive control for power control in islanded DC microgrids," *IECON 2015 - 41st Annual Conference of the IEEE, Industrial Electronics Society*, Yokohama, 2015, pp. 001610-001615.

[26] T. Vandoorn, J. Vasquez, J. D. Kooning, J. M. Guerrero, and L. Vandevelde, "Microgrids: Hierarchical control and an overview of the control and reserve management strategies," *IEEE Industrial Electronics Magazine*, vol. 7, no. 4, pp. 42–55, Dec. 2013.

[27] R. R. Negenborn and J. M. Maestre, "Distributed Model Predictive Control: An Overview and Roadmap of Future Research Opportunities," *IEEE Control Syst.*, vol. 34, pp. 87–97, Aug. 2014.

[28] M. Zangeneh, M. Hamzeh, H. Mokhtari, and H. Karimi, "A new power management control strategy for a MV microgrid with both synchronous generator and inverter-interfaced distributed energy resources," *IEEE*

*Industrial Electronics (ISIE), 2014 IEEE 23rd International Symposium on*, vol., no., pp.2529–2534, June 2014.

[29] W. Tushar, B. Chai, C. Yuen, D. B. Smith, K. L. Wood, Z. Yang, and H. V. Poor, "Three-Party Energy Management With Distributed Energy Resources in Smart Grid," *IEEE Trans. Ind. Electron.*, vol. 62, no. 4, pp. 2487–2498, April 2015.

[30] G. François, and D. Bonvin "Real-Time Optimization: Optimizing the Operation of Energy Systems in the Presence of Uncertainty and Disturbances," *13th International Conference on Sustainable Energy technologies*, Geneva, August 2014.

[31] T. V. Vu, S. Paran, T. El Mezyani, and C. S. Edrington, "Real-time distributed power optimization in the DC microgrids of shipboard power systems," *IEEE Electric Ship Technologies Symposium, ESTS* 2015. pp. 118–122, 2015.

[32] A. N. Venkat, I. A. Hiskens, J. B. Rawlings, S. J. Wright, "Distributed MPC Strategies With Application to Power System Automatic Generation Control," *IEEE Control Syst.*, vol. 16, no. 6, pp. 1192–1206, Nov. 2008.

[33] M. Pipattanasomporn, H. Feroze, and S. Rahman, "Multi-agent systems in a distributed smart grid: Design and implementation," in *Proc. IEEE PES 2009 Power Syst. Conf. Expo. (PSCE)*, Mar. 2009.

[34] C. M. Colson, M. H. Nehrir, and R. W. Gunderson, "Multi-agent microgrid power management," in *Proc. 18th IFAC World Congr.*, Milan, Italy, Sep. 2011.

[35] A. L. Dimeas and N. D. Hatziargyriou, "Operation of a Multiagent System for Microgrid Control," *IEEE Trans. Power Syst.*, vol. 20, no. 3, pp. 1447–1455, Aug. 2005.

[36] C. X. Dou, and B. Liu, "Multi-Agent Based Hierarchical Hybrid Control for Smart Microgrid," *IEEE Trans. Smart Grid*, vol. 4, no. 2, pp. 771-778, June 2013.

[37] Y. Guo, M. Pan, and Y. Fang, "Optimal power management of residential customers in the smart grid," *IEEE Trans. Parallel Distrib. Syst.*, vol. 23, no. 9, pp. 1593–1606, Sep. 2012.

[38] N. Rahbari-Asr, Y. Zhang and M. Y. Chow, "Consensus-based distributed scheduling for cooperative operation of distributed energy resources and storage devices in smart grids," in *IET Generation, Transmission & Distribution*, vol. 10, no. 5, pp. 1268-1277, 4 7 2016.

[39] S Paran, TV Vu, T El Mezyani, CS Edrington, "MPC-based power management in the shipboard power system," *IEEE Electric Ship Technologies Symposium, ESTS 2015*. pp. 118–122, 2015.

[40] L. Yacoubi, K. Al-Haddad, L. A. Dessaint, and F. Fnaiech, "Linear and nonlinear control techniques for a three-phase three-level NPC boost rectifier," *IEEE Trans. Ind. Electron.*, vol. 53, no. 6, pp. 1908–1918, Dec. 2006.

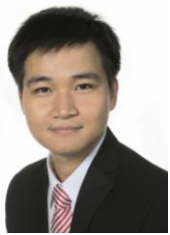
**Tuyen V. Vu** (S'14) was born in Hung Yen, Vietnam. He received his BS in electrical engineering from the Hanoi University of Science Technology, Vietnam in 2012, and his PhD in electrical engineering from the Florida State University in 2016.

He worked as a graduate research assistant for the Florida State University-Center for Advanced Power Systems from 2013 to 2016. Since August 2016, he has been a postdoctoral associate in the Florida State University-Center for Advanced Power Systems. His research interests include modeling, advanced controls, and managements of power electronics, motor drives, microgrids, and power distribution systems.

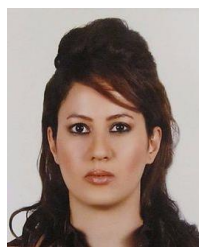
**Sanaz Paran** (S'10) received her B.S. degree in Electrical Engineering from Shiraz University/Iran in 2009, and her MS and PhD in electrical engineering from the Florida State University in 2013 and 2016, respectively.

She has been a research assistant for the Energy Conversion and Integration Thrust at the Center for Advanced Power Systems at Florida State University since 2011. Her research interests include power hardware in the loop applications, applied power electronics, intelligent control of distributed energy resources, and automotive and adaptive energy/power management for AC and DC Microgrids.

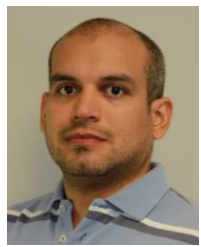
**Fernand Diaz Franco** (S'14) was born in Cali, Colombia. He received the B.S. degree in physics engineering from University of Cauca, Colombia, and the M.Sc. in mechatronics systems from University of Brasilia, Brazil in 2005 and 2008 respectively.

He is currently working towards the Ph.D. degree in electrical engineering at Florida State University.

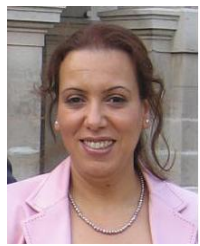
**Touria El-Mezyani** (M'08) received the Ph.D. degree in Control, Computer and Systems Engineering from Université des Sciences et Technologies de Lille, France in 2005.

She is currently Research Faculty with the Center for Advanced, Florida State University, Tallahassee. Her research interests include complex and hybrid systems application to power systems, distributed/decentralized controls and optimization, fault detection, isolation, and diagnosis, control and optimization, multi-agent systems. She is an affiliate of the IEEE Automatic Control Society since 2008.

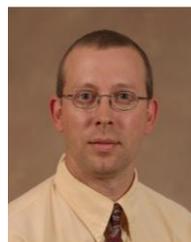
**Chris S. Edrington** (SM' 08) received his BS in engineering from Arkansas State University in 1999 and his MS and PhD in Electrical Engineering from the Missouri Science and Technology (formerly University of Missouri-Rolla) in 2001 and 2004, respectively, where he was both a DoE GAANN and NSF IGERT Fellow.

He currently is a Professor of Electrical and Computer Engineering with the FAMU-FSU College of Engineering and is the lead for the Energy Conversion and Integration thrust for the Florida State University-Center for Advanced Power Systems. His research interests include modeling, simulation, and control of electromechanical drive systems; applied power electronics; distributed control; integration of renewable energy, storage, and pulse power loads.

Dr. Edrington has published over 130 papers (including 2 IEEE Prize Awards), has graduated 25 MS students and 7 PhD students and has 14 patents with and addition patent pending.